\documentclass[journal]{IEEEtran}
\usepackage{cite}

\ifCLASSINFOpdf
  \usepackage[pdftex]{graphicx}
\else
\fi
\usepackage{amsmath}
\usepackage{eurosym}
\usepackage{hyperref}
\usepackage{amssymb}
\usepackage{lipsum, enumitem}
\usepackage{multirow}
\usepackage{xcolor}
\allowdisplaybreaks

\providecommand{\e}[1]{\ensuremath{\times 10^{#1}}}

\usepackage{array}
\newcolumntype{L}[1]{>{\raggedright\let\newline\\\arraybackslash\hspace{0pt}}m{#1}}
\newcolumntype{C}[1]{>{\centering\let\newline\\\arraybackslash\hspace{0pt}}m{#1}}
\newcolumntype{R}[1]{>{\raggedleft\let\newline\\\arraybackslash\hspace{0pt}}m{#1}}

\begin{document}
\NoHyper
%
\title{Second-Order Cone Relaxations of the Optimal Power Flow for Active Distribution Grids: Comparison of Methods}


\author{Lucien Bobo, Andreas Venzke, Spyros Chatzivasileiadis\thanks{L. Bobo, A. Venzke and S. Chatzivasileiadis are with the Department of Electrical Engineering, Technical University of Denmark, Elektrovej 325, Kongens Lyngby, Denmark. E-mail: \{andven, spchatz\}@elektro.dtu.dk\newline\indent This work  is  supported  by  the  EU H2020 project FLEXGRID, funded by the European Commission, Grant Agreement No. 863876.}}
\maketitle


\begin{abstract}
Convex relaxations of the AC Optimal Power Flow (OPF) problem are essential not only for identifying the globally optimal solution but also for enabling the use of OPF formulations in Bilevel Programming and Mathematical Programs with Equilibrium Constraints (MPEC), which are required for solving problems such as the coordination between transmission and distribution system operator (TSO/DSO) or optimal network investment.  Focusing on active distribution grids and radial networks, this paper introduces a framework that collects and compares, for the first time to our knowledge, the performance of the most promising convex OPF formulations for practical applications. Our goal is to establish a solid basis that will inform the selection of the most appropriate algorithm for different applications. This paper (i) introduces a \emph{unified} mathematical and simulation framework, (ii) extends existing methods to  retain exactness in a wider number of cases and (iii) consider reactive power injections. We conduct simulations on the IEEE 34 and 123 radial test feeders with distributed energy resources (DERs), using yearly solar irradiation and load data.

\end{abstract}

%
\IEEEpeerreviewmaketitle
\vspace{-10pt}
{\footnotesize
\section*{Nomenclature}
	\begin{itemize}[leftmargin=11mm]
		\setlength\itemsep{0.1em}
		\item[] \small \textbf{Sets and indices}
		\item[$l$] Index referring to a bus or to a line. Bus $l$ is at the downstream end of line $l$. Bus $0$ is the root node. Bus $1$ is connected to the root node by line $1$.
		\item[$\mathcal{L}$] Set of non-root nodes $l$ or lines $l$ in the distribution grid.
		\item[$\text{up}(l)$] Bus upstream of bus $l$.
		\item[$\mathcal{L}_l$] Set of buses or lines downstream of bus $l$.\\ [-2ex]
	\end{itemize}
	
	\begin{itemize}[leftmargin=11mm]
		\setlength\itemsep{0.1em}
		\item[] \small \textbf{Parameters (physical quantities are in p.u.)}
		\item[$z_l$] Impedance of line $l$.
		\item[$b_l$] Half shunt capacitance of line $l$.
		\item[$v_0$] Squared magnitude of the root node (substation) voltage.
		\item[$v_l^\text{min}$] Squared lower voltage bound at bus $l$.
		\item[$v_l^\text{max}$] Squared upper voltage bound at bus $l$.
		\item[$I_l^\text{max}$] Squared current bound for line $l$.
		\item[$\textbf{G}_{l,m}$] Adjacency matrix of the oriented graph of the network:\\$\textbf{G}_{l,m}=1$ if $l=\text{up}(m)$ and $0$ otherwise.
		\item[$\mathcal{C}_l^\Re,\mathcal{C}_l^\Im$] Active and reactive power withdrawal cost functions, bus $l$.
		\item[$\mathcal{C}_0^\Re,\mathcal{C}_0^\Im$] Active and reactive power import cost functions.
		\item[$S_l$] Feasible region of complex power withdrawal at bus $l$.\\ [-2ex]
	\end{itemize}
	
	\begin{itemize}[leftmargin=11mm]
		\setlength\itemsep{0.1em}
		\item[] \small \textbf{Variables}
		\item[$v_l$] Squared magnitude of the complex voltage at bus $l$.
		\item[$f_{l}$] Squared magnitude of current in central element of line $l$.
		\item[$\bar{v}_{l}, \bar{f}_{l}$] Auxiliary variables.
		\item[$s_l$] Complex power withdrawal at bus $l$: $s_l=p_l+\mathbf{i}q_l$.
		\item[$S_{l}^t$] Complex power flow entering line $l$ from upstream bus $\text{up}(l)$: $S_{l}^t=P_{l}^t+jQ_{l}^t$.
		\item[$S_{l}^b$] Complex power flow entering bus $l$ from line $l$:\item[]$S_{l}^b=P_{l}^b+jQ_{l}^b$.
		\item[$\hat{S}_{l}^t, \bar{S}_{l}^t$] Complex auxiliary variables (real and imag. parts as above).
		\item[$\hat{S}_{l}^b, \bar{S}_{l}^b$] Complex auxiliary variables (real and imag. parts as above).
	\end{itemize}
}	

\section{Introduction}
\label{sec:intro}
Convex relaxations of the AC Optimal Power Flow (OPF) have attracted wide interest over the past years. Several formulations have been proposed with ultimate goal to determine the global optimum of the original OPF problem; for a detailed review see \cite{low14partI,low14partII,molzahn17} and references therein. A convex OPF problem, however, can also be used for a wide range of different applications. Bilevel programs, often proposed for strategic bidding and TSO/DSO coordination, and mathematical programs with equilibrium constraints (MPEC), often employed for optimal investment decisions or optimal operation, \emph{require} the lower level optimization programs to be convex. The electric network constraints represented by a convexified OPF formulation are usually one of the lower problems to be included. Focusing on active distribution grids and radial networks, this paper introduces a framework that collects and compares, for the first time to our knowledge, the performance of the most promising convex OPF formulations for practical applications. 
Our goal is to establish a solid basis that will inform the selection of the most appropriate algorithm for different applications. 

The main challenge of all convex OPF formulations remains their ability to determine an optimal point that is \emph{feasible} to the original problem, i.e. to be exact. The exactness of convex relaxations for meshed networks has been investigated in detail in \cite{venzke19}. In distribution grids (usually radial networks), where the cost of operation is orders of magnitude lower than in transmission grids, determining a feasible operating point is usually more important than identifying the global optimum. Especially in bilevel and MPEC formulations, where a convex OPF formulation is the only option when casting it as a lower level program, a feasible optimal point is almost a requirement. The focus of this paper is, therefore, on methods that can provide guarantees to achieve a feasible optimal solution for \emph{realistic} distribution grids. 

In radial networks, second-order cone relaxation (SOCP) is the preferred formulation as it has been shown to be equivalent to the standard semidefinite relaxation (SDP) \cite{bose12}, while it requires a much lower computational effort than SDP. Although substantial efforts to derive sufficient conditions for exactness of the direct SOCP relaxation have been reported in the literature for radial networks (see \cite{low14partII} for a review), these do not apply in a realistic setting of distribution system operation with high penetration of DERs. 
In particular, the conditions proposed in \cite{gan12} require no upper bounds on voltage levels, while the conditions in \cite{farivar11} require no upper bounds on the withdrawals of active and reactive power. Some authors propose sufficient conditions based on voltage angle differences \cite{lam12,lavaei14}, but these conditions assume fixed voltage levels or no reactive power constraints, and, thus, do not apply to the operation of active distribution networks. 

Ref.~\cite{gan15,huang17,nick17} follow a different approach: they first augment the OPF problem with additional constraints that shrink the feasible space, and subsequently relax the \emph{augmented} problem. This allows them to define milder provable conditions for exactness.
To our knowledge, the methods from \cite{gan15,huang17,nick17} constitute the most promising contributions on \emph{convex} OPF solution methods for radial, single-line networks, that can guarantee a feasible optimal point, even if this might not coincide with the global optimum. Despite each method having compared its performance to the direct SOCP relaxation, so far, no common framework to compare the performance of all three methods has existed and no testing in realistic simulations has taken place to examine the extent of their possible practical application. This is the first paper, to our knowledge, that collects all three methods under a common mathematical and simulation framework, and provides a thorough comparison of their performance. The contributions of this paper are the following: 
\begin{enumerate}
    \item We collect the formulations, assumptions, and conditions for exactness of \cite{gan15,huang17,nick17} under a \emph{unified mathematical framework}. We assess their exactness and suboptimality for wide range of realistic conditions, generating hourly load and PV profiles for an entire year on the IEEE34 and IEEE123 test feeders. We provide the feeder data, with the required modifications that facilitated the comparison, in an on-line appendix for future use by the interested reader \cite{appendix}. 
    \item For practical applications that violate the stipulated sufficient conditions, we introduce a mild penalty on current magnitudes in the objective function to achieve exactness without incurring substantial sub-optimality.  With the help of the penalty factor, we find that the methods \cite{gan15,huang17,nick17} can be successfully applied to a substantially wider range than the theoretically sufficient conditions suggest.
   
    \item We extend the original range of applications of \cite{gan15,huang17,nick17} beyond active power, by introducing reactive power control, and assess the exactness of the examined methods. 
    \item We specify the use of a numerical threshold to assess exactness. In our simulation results, it appears that when all residuals of the relaxed constraints are below $1\e{-2}$ p.u, the relaxation is exact. 
   	\item We formulate recommendations on the best use of the reviewed methods for practical applications, and identify gaps to be filled to extract milder sufficient conditions for exactness.
\end{enumerate}
The rest of this paper is structured as follows. Section \ref{sec:methods} formulates an OPF problem which generalises the models in \cite{gan15,huang17,nick17} and compares assumptions and augmentations of the three methods. Section \ref{sec:simu} explains our methodology for evaluating exactness, introduces the current penalty term, and presents the test case setup. Section \ref{sec:results} presents simulation results and main insights. Finally, Section \ref{sec:ccl} formulates recommendations and concludes.

\section{Optimal Power Flow Formulation}
\label{sec:methods}
As already mentioned, when it comes to using optimal power flow as a decision support tool in real applications, it is crucial that the OPF determines feasible operating points. The methods in \cite{gan15,huang17,nick17} seek to formulate a convex optimization problem whose solution is guaranteed to be inside the feasible region of the original OPF. 
The methods differ in the notation they use, the simplifications they make, the augmentation of the OPF problem they propose, and the conditions for exactness they define. Some important differences between the methods are summarised in Table \ref{tab:methods}.
\begin{table}
\caption{Assumptions and conditions for exactness in \cite{gan15,huang17,nick17}}
\label{tab:methods}
\begin{center}
\vspace{-10pt}
\begin{tabular}{|c|c|c|c|}
\hline 
• & \cite{gan15} & \cite{huang17} & \cite{nick17} \\ 
\hline 
models line shunt capacitance & no & no & yes \\ 
\hline 
considers line current limits & no & yes & yes \\ 
\hline 
considers apparent power flow limits & no & yes & no \\ 
\hline 
same voltage bounds at each bus & no & no & yes \\ 
\hline 
conditions on objective function & (a)(c) & (a)(d) & (b)(c) \\
\hline
conditions on network parameters & yes & no & yes \\
\hline
\end{tabular}\\
\end{center}
\vspace{-3pt}
\begin{itemize}
\item[(a)] Objective sums cost functions for active power import at the root node and for active power withdrawal at non-root nodes.
\item[(b)] Objective sums cost functions for active power import at the root node and for active and reactive power withdrawal at non-root nodes.
\item[(c)] Cost function for active power import is strictly increasing.
\item[(d)] Cost function for active power import is non-decreasing.
\vspace{-15pt}
\end{itemize}
\end{table}
All methods formulate the OPF problem using the branch flow model, and all assume the root node voltage to be a fixed parameter. In the following, Section \ref{sec:opf} presents an OPF model which unifies the formulations from \cite{gan15,huang17,nick17}. Section \ref{sec:soc} introduces its SOCP relaxation, and \ref{sec:aug} describes the augmentations proposed in \cite{gan15,huang17,nick17}.
\subsection{Non-Convex Optimal Power Flow}
\label{sec:opf}
A radial distribution network consists of lines and buses with one root node connected to the higher-voltage grid. Using a branch flow model, with lines represented by a $\pi$-model similar to that of \cite{nick17}, let \eqref{OPFOF}-\eqref{OPF11} describe the non-convex OPF problem. 
The objective \eqref{OPFOF} accounts for costs associated with the active and reactive power withdrawal at non-root node $l$, as well as the import of active and reactive power from the main grid through the root node; Equations \eqref{OPF1}-\eqref{OPF4} model the load flow using the branch flow model; Equations \eqref{OPF5}-\eqref{OPF11} impose bounds on the bus voltages, line currents and apparent power flows, and enforce the complex power withdrawal at each non-root node to remain with a given feasible region $\mathcal{S}_l\subset\mathbb{C}$. The region $\mathcal{S}_l$ is bounded from below, i.e., there exist $s^\text{min}_l$ such that $\mathcal{S}_l\subseteq\{s_l\in\mathbb{C}~|~s_l\geq s^\text{min}_l\}$, for all $l\in\mathcal{L}$. Also, let $s_l=p_l+jq_l$, $S_{l}^t=P_{l}^t+jQ_{l}^t$ and $S_{l}^b=P_{l}^b+jQ_{l}^b$.
\begin{align}
\label{OPFOF}
& \min_{\substack{s_l,v_l,f_l\\S_l^b,S_l^t}}~\sum_{l\in\mathcal{L}}\left( \mathcal{C}^\Re_l\left(p_l\right) + \mathcal{C}^\Im_l\left(q_l\right) \right)+ \mathcal{C}^\Re_0(P_1^t) + \mathcal{C}^\Im_0(Q_1^t),\\[-1ex]
\label{OPF1}
& \quad\text{s.t.}~~\bigg\{ S_l^b=s_l+\sum_{m\in\mathcal{L}}(\textbf{G}_{l,m}S_m^t),\\
\label{OPF2}
& \quad S_l^t=s_l+\sum_{m\in\mathcal{L}}(\textbf{G}_{l,m}S_m^t) +z_l f_l-j(v_{\text{up}(l)}+v_l)b_l,\\
\label{OPF3}
& \quad v_l = v_{\text{up}(l)} - 2\Re\left(z_l^\ast(S_l^t+j v_{\text{up}(l)}b_l)\right)+|z_l|^2f_l,\\
\label{OPF4}
& \quad f_l=\frac{|S_l^t+j v_{\text{up}(l)} b_l|^2}{v_{\text{up}(l)}},\\
\label{OPF5}
& \quad v_l^\text{min}\leq v_l\leq v_l^\text{max},\\
\label{OPF6}
& \quad |S_l^b|^2 \leq I_l^\text{max} v_l,~|S_l^t|^2 \leq I_l^\text{max} v_{\text{up}(l)},\\
\label{OPF9}
& \quad |S_l^b| \leq S_l^\text{max},~|S_l^t| \leq S_l^\text{max},\\
& \quad s_l \in \mathcal{S}_l \bigg\}, \forall~l\in\mathcal{L} .
\label{OPF11}
\vspace{-0.2cm}
\end{align}
Under the respective restrictions listed in Table \ref{tab:methods}, the OPF formulations from \cite{gan15,huang17,nick17} are equivalent to \eqref{OPFOF}-\eqref{OPF11}.
\vspace{-0.2cm}
\subsection{SOCP relaxation}
\label{sec:soc}
\begin{figure}
\centering
\includegraphics[width=3.3in]{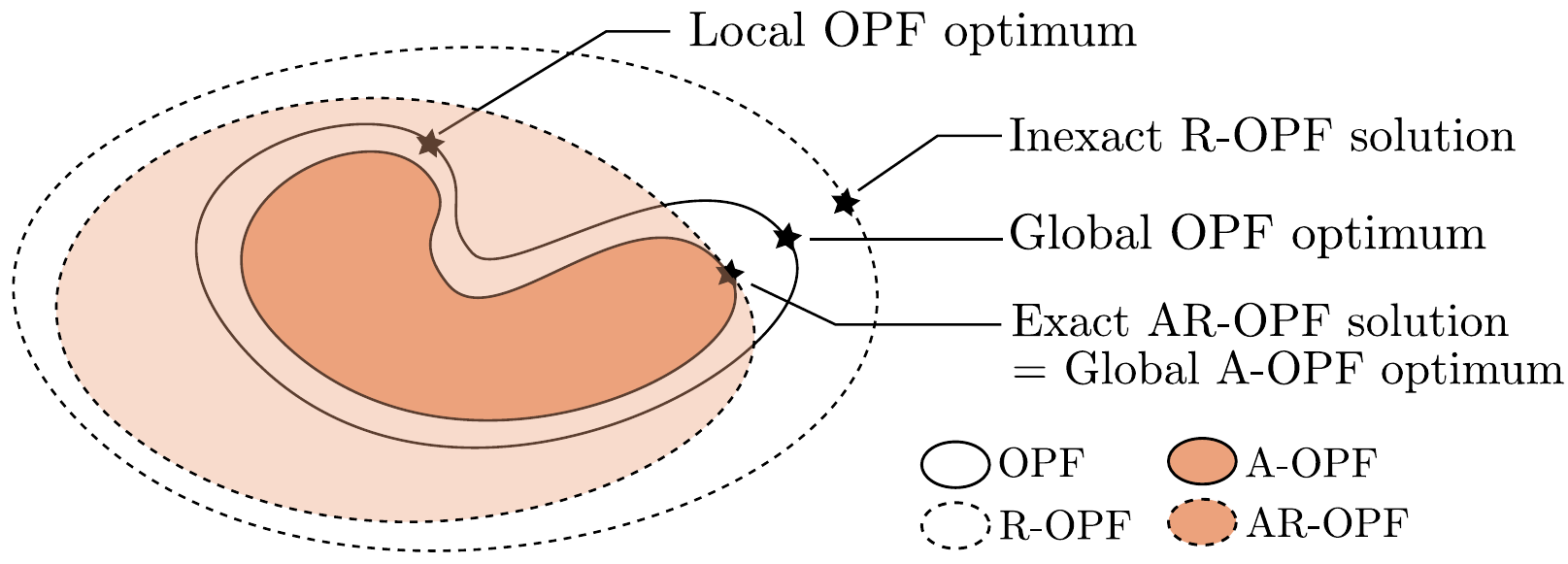}
\vspace{-5pt}
\caption{Venn diagram of a case where the relaxation (R-OPF) is inexact and the augmented relaxation (AR-OPF) is exact (set sizes are not drawn to scale).}
\label{fig:venn}
\vspace{-10pt}
\end{figure}

The OPF problem \eqref{OPFOF}-\eqref{OPF11} is non-convex, but can be relaxed into an SOCP problem by substituting constraint (\ref{OPF4}) with the following second-order cone constraint:
\begin{equation}
\label{eq:SOC}
v_{\text{up}(l)} f_l \geq |S_l^t+j v_{\text{up}(l)} b_l|^2, ~\forall l \in\mathcal{L},
\end{equation}
The relaxed OPF problem (R-OPF) writes:
\begin{align}
\label{SOCPbaseOF}
&\min_{\substack{s_l,v_l,f_l\\S_l^b,S_l^t}}~\sum_{l\in\mathcal{L}}\left( \mathcal{C}^\Re_l\left(p_l\right) + \mathcal{C}^\Im_l\left(q_l\right) \right)+ \mathcal{C}^\Re_0(P_1^t) + \mathcal{C}^\Im_0(Q_1^t),\\[-1ex]
\label{SOCPbaseEQ}
&\quad\text{s.t.}~~\text{\eqref{OPF1}-\eqref{OPF3}, \eqref{eq:SOC}, \eqref{OPF5}-\eqref{OPF11}}, ~\forall~l\in\mathcal{L}.
\end{align}
As discussed in Section \ref{sec:intro}, this direct relaxation of the OPF problem is only proven to be exact under sufficient conditions that are difficult to meet in practice for active distribution networks \cite{low14partII,gan12,farivar11,lam12,lavaei14}. The authors in \cite{gan15,huang17,nick17} propose instead to solve \textit{augmented relaxations} of the OPF problem (referred to as \textit{AR-OPF}) which retain exactness under milder conditions. These are based on augmented versions of the non-convex OPF problem (referred to as \textit{A-OPF}), which expand \eqref{OPF1}-\eqref{OPF11} with additional constraints (see Figure \ref{fig:venn} for a schematic representation of the different problems). The relaxation step is similar to R-OPF, i.e. the second-order cone constraint \eqref{eq:SOC} substitutes the non-convex constraint \eqref{OPF4}. 
Note that, with OF$^\star$ denoting the \emph{globally} optimal objective function of a given problem,
the following relation holds when AR-OPF is exact:
\begin{gather}
\text{OF}^\star(\text{OPF}) \leq \text{OF}^\star(\text{AR-OPF}) = \text{OF}^\star(\text{A-OPF})
\end{gather}
\vspace{-0.9cm}
\subsection{Augmentations of the OPF}\label{sec:aug}
Gan et al. \cite{gan15} propose to augment the OPF problem and its relaxation with (\ref{SOCP_Gan1})-(\ref{SOCP_Gan3}), $\forall l \in\mathcal{L}$. The authors prove that their AR-OPF is exact under a sufficient ex-ante condition computed from the network parameters (see \cite{gan15}) if the active power import cost function $\mathcal{C}_0^\Re$ is strictly increasing. Their OPF formulation omits line shunts, as well as line current limits and apparent power flow limits (see Table \ref{tab:methods}).
\begin{align}
\label{SOCP_Gan1}
& \hat{S}_l^b=s_l+\sum_{m\in\mathcal{L}}(\textbf{G}_{l,m}\hat{S}_m^b),\\
\label{SOCP_Gan2}
& \bar{v}_l=\bar{v}_{\text{up}(l)}-2\Re\left(z_l^\ast\hat{S}^b_l\right),\\
& \bar{v}_l \leq v^\text{max}_l.
\label{SOCP_Gan3}
\end{align}




Huang et al. \cite{huang17} further propose to augment the OPF problem and its relaxation with (\ref{SOCP_Gan1})-(\ref{SOCP_Gan3}) and \eqref{SOCP_Huang1}. The authors prove that the withdrawal quantities $s_l$ at an optimal solution of their AR-OPF lead to a feasible load flow (a property implied by exactness, but not implying exactness as defined in this paper, see Section \ref{sec:methex}) if the active power import cost function $\mathcal{C}_0^\Re$ is non-decreasing. Their OPF formulation includes line current and apparent power flow limits, and omits line shunts.
\begin{align}
\label{SOCP_Huang1}
& \Re\left(z_m^\ast\hat{S}^b_l\right)\geq 0,&\forall l \in\mathcal{L}, m \in \mathcal{L}_l.
\end{align}



Nick et al. \cite{nick17} propose to augment the OPF problem and its relaxation with (\ref{SOCP_Nick1})-(\ref{SOCP_Nick11}), $\forall l \in\mathcal{L}$, where $P_l^\text{max}$ and $Q_l^\text{max}$ are parameters chosen so that they do not restrict the feasible space of the A-OPF. The authors prove that their AR-OPF is exact under sufficient ex-ante conditions computed from the network parameters (see \cite{nick17}) if the active power import cost function $\mathcal{C}_0^\Re$ is strictly increasing. Their formulation includes line shunts as well as line current limits, omits apparent power flow limits, and assumes $v^\text{max}_l$ have the same value for all $l$. When line shunts and line current limits are ignored and $v^\text{max}_l$ have the same value for all $l$, the augmentations of \cite{gan15} and \cite{nick17} are equivalent.
\begin{align}
\label{SOCP_Nick1}
& \hat{S}_l^t=s_l+\sum_{m\in\mathcal{L}}(\textbf{G}_{l,m}\hat{S}_m^t)-j(\bar{v}_{\text{up}(l)}+\bar{v}_l)b_l,\\
\label{SOCP_Nick2}
& \bar{v}_l = \bar{v}_{\text{up}(l)} - 2\Re\left(z_l^\ast(\hat{S}_l^t+j \bar{v}_{\text{up}(l)}b_l)\right),\\
\label{SOCP_Nick3}
& \bar{S_l^t}=s_l+\sum_{m\in\mathcal{L}}(\textbf{G}_{l,m}\bar{S}_m^t)+z_l\bar{f}_l\\
& \bar{f}_l v_l \geq \text{max}\left\lbrace(\hat{P}_l^b)^2,(\bar{P}_l^b)^2\right\rbrace&\nonumber\\
\label{SOCP_Nick4}
&~~~~~~~~~~ + \text{max}\left\lbrace(\hat{Q}_l^b-\bar{v}_lb_l)^2,(\bar{Q}_l^b-v_lb_l)^2\right\rbrace,\\
& \bar{f}_l v_{\text{up}(l)} \geq \text{max}\left\lbrace(\hat{P}_l^t)^2,(\bar{P}_l^t)^2\right\rbrace&\nonumber\\
\label{SOCP_Nick5}
&~~~~~~~~~~ + \text{max}\left\lbrace(\hat{Q}_l^t-\bar{v}_{\text{up}(l)}b_l)^2, (\bar{Q}_l^t-v_{\text{up}(l)}b_l)^2\right\rbrace,\\
\label{SOCP_Nick6}
& \bar{S}_l^b=s_l+\sum_{m\in\mathcal{L}}(\textbf{G}_{l,m}\bar{S}_m^t),~~ \hat{S}_l^b=s_l+\sum_{m\in\mathcal{L}}(\textbf{G}_{l,m}\hat{S}_m^t),\\
\label{SOCP_Nick7}
& \bar{v}_l \leq v^\text{max}_l,\\
\label{SOCP_Nick9}
& \text{max}\left\lbrace\hat{P}_l^b,\bar{P}_l^b\right\rbrace^2 + \text{max}\left\lbrace\hat{Q}_l^b,\bar{Q}_l^b\right\rbrace^2 \leq v_l I_l^\text{max},\\
\label{SOCP_Nick10}
& \text{max}\left\lbrace\hat{P}_l^t,\bar{P}_l^t\right\rbrace^2+ \text{max}\left\lbrace\hat{Q}_l^t,\bar{Q}_l^t\right\rbrace^2 \leq v_l I_l^\text{max},\\
\label{SOCP_Nick11}
& P_l^t \leq \bar{P}_l^t \leq P_l^\text{max},~ Q_l^t \leq \bar{Q}_l^t \leq Q_l^\text{max}.
\end{align}

\section{Methodology}\label{sec:simu}
\subsection{Evaluating Exactness}\label{sec:methex}
A relaxation (R-OPF or AR-OPF) is said to be \textit{exact} if its optimal solutions are feasible to the non-relaxed problem, i.e. satisfy the power flow equations and all inequality constraints. In this case, these solutions are globally optimal solutions to the non-relaxed problem (resp. OPF or A-OPF). For the SOCP relaxations studied in this paper, the relaxed problem is exact if the inequalities \eqref{eq:SOC} are binding at optimality, i.e. the residuals of the relaxed constraints are zero. This definition is in line with the use of the term in \cite{gan15} and \cite{nick17}, and is also used in this paper (see below for the definition used in \cite{huang17}). 
Due to numerical accuracy of existing SOCP solvers, however, residuals are unlikely to be exactly zero in simulation results. In our numerical analyses, we consider a relaxation to be exact when all residuals of the relaxed constraints are relatively close to zero -- more specifically below $1\e{-2}$ p.u., which appeared to work well as a threshold in our simulations. 

To confirm that the optimal solutions to an exact relaxed problem can be used in practice, we compute a load flow based on the active and reactive power setpoints $s_l$ prescribed by the relaxed problem (R-OPF or AR-OPF). We then evaluate whether the voltage and current levels determined by the load flow solution are exceeding their bounds. We refer to the \textit{voltage bound violation} of a load flow solution as the largest positive difference between voltage levels and their upper bounds, in nominal values, i.e. $\max\{0,(\sqrt{v^{\textsc{lf}}_l}-\sqrt{ v^{\text{max}}_l})_{l\in\mathcal{L}}\}$, where $v^{\textsc{lf}}_l$ is the squared voltage magnitude at bus $l$ in the load flow solution. \textit{Current bound violation} is defined analogously. With numerical accuracy in mind, we consider that a solution can be used in practice if both voltage and current bound violation are below $1\e{-2}$ p.u.. Note that the authors in \cite{huang17} refer to this property as exactness, considering, in other words, the R-OPF/AR-OPF solution as a warm-start to a load flow computation. If the load flow computation results to voltage and current levels that are within bounds, then they consider the relaxation exact.

\subsection{Penalty term on current magnitudes}
In inexact cases, we measure the \textit{relaxation gap} as the difference between the optimal objective value of the relaxed problem and the optimal objective value of the original problem, i.e. $\text{OF}^\star(\text{AR-OPF}) - \text{OF}^\star(\text{A-OPF})$ for an augmented relaxation, or $\text{OF}^\star(\text{R-OPF}) - \text{OF}^\star(\text{OPF})$ for the original OPF. In our simulations, we observed cases that, despite achieving a zero relaxation gap, remain inexact. As we will also see later in this paper, this occurs because of a flat objective function, where several operating points achieve the same objective function value, and the optimisation is unable to determine a feasible one (see also \cite{venzke19}).

Penalisation methods are commonly used to tighten inexact relaxations (see for example \cite{madani2014convex, zohrizadeh2018penalized}), especially in cases where solution feasibility is more important than solution optimality. There are however configurations where penalising the objective function can yield feasible solutions without having to compromise with optimality. In this paper, for cases where the solution of a relaxed problem is inexact with a zero relaxation gap (i.e. there exists multiple optimal solutions, some of which are inexact), we propose to add a mild penalty term on the squared current magnitudes $f_l$ with weight $\epsilon$ to the objective \eqref{SOCPbaseEQ} to drive the optimisation towards an exact solution:
\begin{align}
    \sum_{l\in\mathcal{L}}\left( \mathcal{C}^\Re_l\left(p_l\right) + \mathcal{C}^\Im_l\left(q_l\right) \right)+ \mathcal{C}^\Re_0(P_1^t) + \mathcal{C}^\Im_0(Q_1^t) + \epsilon f_l 
\end{align} Numerical examples in Section~\ref{res:4} show that a small penalty weight allows to recover an exact solution with negligible sub-optimality.
\subsection{Test Case Setup}
\label{sec:testcases}

We evaluate the solution methods using two test feeders, IEEE34 and IEEE123 \cite{feeders}, to which we add distributed generation in the form of PV stations. IEEE34 is a rather small network with long line segments, while IEEE123 is a larger network with both overhead and underground line segments. These test networks are originally unbalanced three-phase systems and include transformers, breakers, capacitor banks, voltage regulators, as well as a mix of spot and distributed loads. As the OPF model \eqref{OPFOF}-\eqref{OPF11} presented in the previous section does not account explicitly for these components, we make a series of adjustments outlined below.

In this study, the networks are assumed to be balanced, and all loads are assumed to be spot PQ loads. We equally distribute unbalanced loads over the three phases, and the distributed loads over the neighbouring buses. We model the feeders as single-line systems using the positive-sequence impedance of the lines for one of their phases, and multiply the loads on single-phase lines by a factor of three to adequately account for voltage drops. For both networks, we use a three-phase power base value of 1 MVA, a phase-to-phase voltage base value of 4.16 kV for IEEE123 and the low-voltage part of IEEE34, and of 24.9 kV for the rest of IEEE34. Transformers are modelled as a line with series resistance and inductance, following the IEEE specifications. Capacitor banks are modelled as fixed reactive power injections. Voltage regulators and breakers are omitted, and the buses upstream and downstream of each voltage regulator in IEEE34 are merged. Data for the feeder models used in our simulations is provided in an on-line appendix \cite{appendix}.

For both networks, we run simulations for all hours of a year-long load profile, where the loads fluctuate between 15\% and 100\% of their default value.\footnote{Data from \url{https://openei.org/}, see on-line appendix \cite{appendix}.} We also add PV stations at all buses with loads connected to them, with a total installed capacity equal to 250\% of the total peak load in each network. The active power output of each PV station is a positive control variable in the OPF, bounded by the hourly available active power production following a year-long profile.\footnote{Data from \url{https://www.nrel.gov/}, see on-line appendix \cite{appendix}.} The reactive power output of each PV station is a free control variable in the OPF, and the apparent power output of each PV station is bounded by a nameplate capacity set to 110\% of its peak active power capacity \cite{turitsyn11}. Unless otherwise specified, the objective function in the simulations is to minimise the import of active power through the root node. In order to compare the different methods based on appropriate modelling assumptions, we define three different simulation configurations: \textit{NS-NC} omits line shunts and current bounds; \textit{NS-C} omits line shunts but includes current bounds; \textit{S-C} includes line shunts and current bounds. 
Unless otherwise specified, we use a fixed root node voltage of 1 p.u., voltage bounds of 0.9 and 1.1 p.u. at all buses (i.e. $v_l^{min}=0.81$ p.u., $v_l^{max}=1.21$ p.u.), and current bounds of 4 p.u. on all lines (i.e. $I_l^{max}=16$ p.u.). We omit apparent power flow limits. 
Optimisation problems are implemented in Matlab using YALMIP \cite{yalmip}. Non-convex OPF and load flows are solved using IPOPT \cite{ipopt} and SOCP are solved using MOSEK \cite{mosek}.


\section{Results}
\label{sec:results}

This section presents simulation results. First, we investigate the feasible space reduction due to the different augmentations and the resulting suboptimality. Second, we assess the conservativeness of the proposed sufficient conditions for exactness in practical applications. Third, we assess the applicability of the relaxations for reactive power control. Finally, we show results for tightening relaxations with a zero relaxation gap with the proposed current penalty term.

\subsection{Shrinkage of the feasible space due to augmentation}
\label{res:1}

As illustrated in Fig. \ref{fig:venn}, the methods from \cite{gan15,huang17,nick17} augment the original OPF problem to alter its feasible region (A-OPF). Through that, their goal is to bring the global optimum of those problems \emph{inside} the OPF original feasible region. In order to determine the A-OPF global optimum, they solve its SOC relaxation (AR-OPF) and derive sufficient conditions to guarantee the AR-OPF exactness. The resulting A-OPF global optimum may be, but is not necessarily, the global optimum to the OPF problem. In this Section, we investigate how suboptimal the A-OPF global optimum is compared with the solution obtained by the non-convex OPF (we know that the non-convex OPF solution is a global optimum when the R-OPF optimal solution has the same objective value). We also identify cases where the feasible region of the A-OPF is empty, i.e. where the augmentation cannot be used to recover a feasible solution to the OPF. Results are summarised in Table \ref{tab:VA}.
\subsubsection{Augmentation from Gan et al. \cite{gan15}}
In this augmentation, described by \eqref{SOCP_Gan1}-\eqref{SOCP_Gan3}, constraint \eqref{SOCP_Gan3} applies the upper voltage bounds $v^\text{max}_l$ to auxiliary variables $\bar v_l$ which are, by construction, greater than or equal to $v_l$ (see proof in \cite{gan15}). This implies that OPF solutions with tight voltage bounds may not be feasible to the A-OPF from Gan et al. \cite{gan15}, as the bound in \eqref{SOCP_Gan3} may be violated when the voltage upper bound in \eqref{OPF5} is binding for one or several buses. In such cases, optimal voltage levels in A-OPF would be lower than in OPF, and exact solutions from the AR-OPF may thus be sub-optimal as compared to the global OPF optimum.

With the NS-NC-34 test case, the AR-OPF from Gan et al. \cite{gan15} finds a global optimum to the OPF in some hours, but yields sub-optimal solutions in others (additional import of active power in those hours is up to 12.8\% of the grid's peak load as compared to the solution of the non-convex OPF). Sub-optimality can be high in the IEEE34 test case because it features heavy voltage drops. Throughout the simulated year, however, this test case only leads to an average additional import of 0.87\% of the grid's peak load as compared to the solution of the non-convex OPF, as upper voltage bounds are tight at optimality only in hours where the available active power from the distributed PV stations is high. Sub-optimality is negligible in the NS-NC-123 test case, where the network is less stressed. Results are summarised in Table \ref{tab:VA}.

We conclude that the method from \cite{gan15} can yield non-negligible sub-optimality in specific cases, but these may be rare enough for the method's applicability not to be compromised.

\subsubsection{Augmentation from Nick et al. \cite{nick17}}
In contrast with \cite{gan15}, the OPF formulation from Nick et al. \cite{nick17} includes line shunts and current bounds (see Table \ref{tab:methods}). The augmentation from \cite{nick17} is described by Equations \eqref{SOCP_Nick1}-\eqref{SOCP_Nick11}. If we neglect line shunts and current bounds, then \cite{nick17} becomes equivalent to  \cite{gan15}.

Constraints \eqref{SOCP_Nick7}-\eqref{SOCP_Nick10} apply both voltage and current bounds to terms that are, by construction, greater than or equal to the actual voltage and current levels (see proof in \cite{nick17}). 
This implies that OPF solutions with tight voltage or current bounds may not be feasible to the A-OPF, as constraints \eqref{SOCP_Nick7}-\eqref{SOCP_Nick10} may be violated when the upper bounds in \eqref{OPF5}-\eqref{OPF6} are binding for one or several buses or lines. In such cases, optimal voltage or current levels in A-OPF would be lower than in OPF, and exact solutions from the AR-OPF may thus be sub-optimal as compared to the global OPF optimum. With the S-C-34 test case, sub-optimality varies between 0 and 12.7\% of the grid's peak load, with an average of 0.84\% (results resemble those obtained with \cite{gan15} on NS-NC-34 as the effect of line shunts and chosen current bounds on the IEEE34 feeder are mild). With the S-C-123 test case, sub-optimality varies between 0 and 3.8\% of the grid's peak load, with an average of 0.6\% (results differ from those obtained with \cite{gan15} on NS-NC-123 as line shunts are substantial in this feeder, and chosen current bounds constrain the optimal power flows in multiple hours). Results are summarised in Table \ref{tab:VA}.

We conclude that the method from Nick et al. \cite{nick17} can lead to non-negligible sub-optimality in specific conditions, but these conditions may be rare enough for the method's applicability not to be compromised.
\begin{table*}
\caption{Effect of augmentations on feasible space. Peak and average sub-optimality are expressed in \% of the grid's peak load. Infeasibility is expressed in \% of infeasible hours in the simulation year.}
\vspace{-10pt}
\label{tab:VA}
\begin{center}
\begin{tabular}{l}
\begin{tabular}{|c|c|c|c|c|c|c|c|c|c|c|c|c|c|c|}
\hline test grid & \multicolumn{7}{c|}{IEEE34} & \multicolumn{7}{c|}{IEEE123} \\ 
\hline
configuration & NS-NC & NS-NC & \hspace{-5pt}NS-NC$^\ast$\hspace{-5pt} & NS-C & \hspace{-5pt}NS-C$^\ast$\hspace{-5pt} & \hspace{-5pt}NS-C$^\ast$\hspace{-5pt} & S-C & NS-NC & NS-NC & \hspace{-5pt}NS-NC$^\ast$\hspace{-5pt} & NS-C & \hspace{-5pt}NS-C$^\ast$\hspace{-5pt} & \hspace{-5pt}NS-C$^\ast$\hspace{-5pt} & S-C \\
\hline AR-OPF & \hspace{-5pt}\cite{gan15},\cite{nick17}\hspace{-5pt} & \hspace{-5pt}\cite{huang17}\hspace{-5pt} & \hspace{-5pt}\cite{huang17}\hspace{-5pt} & \hspace{-5pt}\cite{huang17}\hspace{-5pt} & \hspace{-5pt}\cite{huang17}\hspace{-5pt} & \hspace{-5pt}\cite{nick17}\hspace{-5pt} & \hspace{-5pt}\cite{nick17}\hspace{-5pt} & \hspace{-5pt}\cite{gan15},\cite{nick17}\hspace{-5pt} & \hspace{-5pt}\cite{huang17}\hspace{-5pt} & \hspace{-5pt}\cite{huang17}\hspace{-5pt} & \hspace{-5pt}\cite{huang17}\hspace{-5pt} & \hspace{-5pt}\cite{huang17}\hspace{-5pt} & \hspace{-5pt}\cite{nick17}\hspace{-5pt} & \hspace{-5pt}\cite{nick17}\hspace{-5pt} \\
\hline \hspace{-5pt}peak sub-opt.\hspace{-5pt} & 12.8\% & 106\% & 124\% & \textcolor{black}{104\%} & \textcolor{black}{122\%} & \textcolor{black}{16.1\%} & 12.7\% & 0.26\% & - & 145\% & \textcolor{black}{-} & \textcolor{black}{68\%} & \textcolor{black}{3.7\%} & 3.8\% \\
\hline \hspace{-5pt}avg. sub-opt.\hspace{-5pt} & 0.87\% & \hspace{-5pt}9.4\%$^{\ast\ast}$\hspace{-5pt} & \hspace{-5pt}24\%$^{\ast\ast}$\hspace{-5pt} & \textcolor{black}{9.9\%$^{\ast\ast}$} & \hspace{-5pt}\textcolor{black}{24\%$^{\ast\ast}$}\hspace{-5pt} & \hspace{-5pt}\textcolor{black}{0.84\%}\hspace{-5pt} & \hspace{-5pt}0.84\%\hspace{-5pt} & \hspace{-5pt}0.006\%\hspace{-5pt} & - & 21\% & \textcolor{black}{-} & \textcolor{black}{19\%} & \textcolor{black}{0.6\%} & 0.6\% \\
\hline infeasibility & 0\% & 90\% & 13\% & \textcolor{black}{90\%} & \textcolor{black}{13\%} & \textcolor{black}{0\%} & 0\% & 0\% & 100\% & 0\% & \textcolor{black}{100\%} & \textcolor{black}{0\%} & \textcolor{black}{0\%} & 0\% \\
\hline
\end{tabular}
\vspace{2pt}\\
$^\ast$ Capacitor banks with variable capacitance.\\
$^{\ast\ast}$ Computed only for feasible hours.
\end{tabular}
\end{center}
\vspace{-10pt}
\end{table*}

\subsubsection{Augmentation from Huang et al. \cite{huang17}}

\begin{figure}[]
\centering
\includegraphics[width=3.2in]{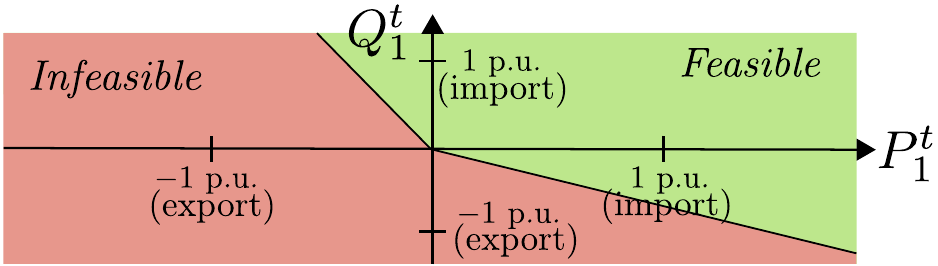}
\vspace{-5pt}
\caption{Feasible space (green area) of power flows allowed by constraint \eqref{SOCP_Huang1}, for the line connected to the root node in IEEE34 test feeder. As the Figure shows, not only simultaneous reverse flows are infeasible, but also any reverse flow (resp. active or reactive) of larger magnitude.}
\label{fig:reverse}
\vspace{-5pt}
\end{figure}

In contrast with \cite{nick17}, the A-OPF from \cite{huang17} directly applies current bounds to current variables $f_l$, meaning that solutions where current levels are close to their limits are not cut off in the augmentation. On the other hand, this augmentation, described by Equations \eqref{SOCP_Gan1}-\eqref{SOCP_Huang1}, constrains reverse and reactive power flows on non-leaf lines (as an illustration, see Figure \ref{fig:reverse} for the power flows allowed by constraint \eqref{SOCP_Huang1} on the line connected to the root node in the IEEE34 feeder). In contrast with the physical interpretation provided by the authors in \cite{huang17}, this augmentation not only prevents simultaneous reverse active and reactive power flows, but also cuts off a substantial share of the feasible space for reverse active or reactive power flows. 
As argued by the authors in \cite{huang17}, this makes their method best suitable for applications with no or light reverse power flows: in a network with strong penetration of non-curtailable DERs, the AR-OPF from \cite{huang17} may be infeasible; when DERs are curtailable, however, the AR-OPF may provide exact, but largely sub-optimal solutions due to shrinkage of the feasible space, as shown in Table~\ref{tab:VA}. 
Our test cases further show that the existence of capacitor banks can be sufficient to render the method sub-optimal or infeasible. While the active and reactive power production from DERs is curtailable in our test cases, the AR-OPF from \cite{huang17} is feasible only for 10\% of hours in the NS-C-34 test case, and never with the NS-C-123 test case. Instead, when allowing the capacitance of the capacitor banks to vary between zero and their original capacitance value, the AR-OPF from \cite{huang17} is feasible in 87\% of hours with NS-C-34, and in all hours with NS-C-123. We conclude that the method from Huang et al. \cite{huang17} is not suited for networks where either DERs or capacitor banks lead to necessary or valuable strong reverse flows. 

\subsection{Exactness of the relaxations}
\label{res:2}

The required conditions so that the direct SOCP relaxation (R-OPF) and the augmented relaxations (AR-OPF) \cite{gan15,huang17,nick17} can yield an exact solution can be in practice quite constraining. In this Section, we evaluate their conservativeness by testing the ability of the four methods to yield exact solutions in practical applications which do not satisfy the proposed sufficient conditions.

\subsubsection{Direct relaxation}
Although a direct SOCP relaxation (R-OPF) has been proven to be exact under some conditions \cite{gan12, farivar11}, the required mathematical conditions for exactness are unpractical, requiring notably the absence of binding upper bounds for voltage levels \cite{gan12} or active and reactive power withdrawals \cite{farivar11}. 
Before moving on with the augmented solution methods (AR-OPF) from \cite{gan15,huang17,nick17}, in this paragraph we examine the exactness of the R-OPF when the proposed sufficient conditions do not hold.

\begin{figure}[]
\centering
\includegraphics[width=3.5in]{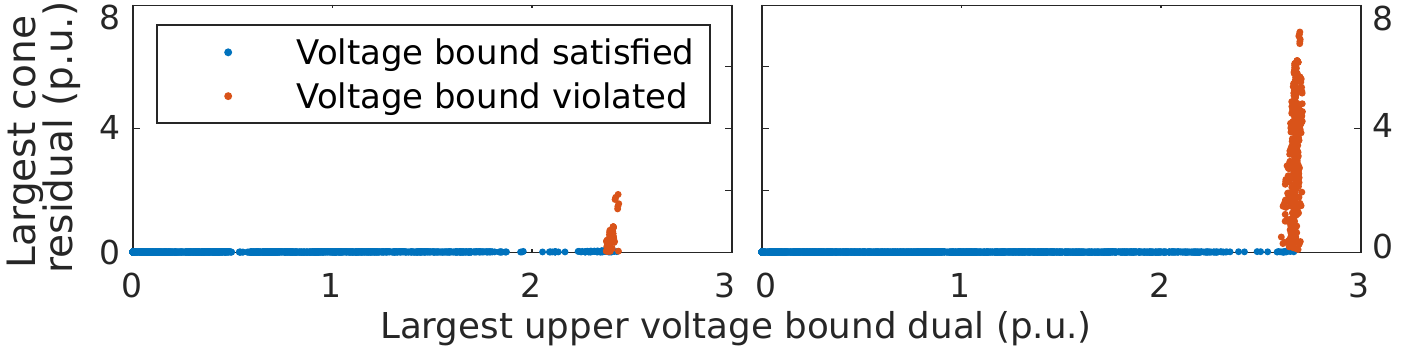}
\vspace{-15pt}
\caption{Largest cone residual in R-OPF solution as a function of the largest dual variable of upper voltage bound. Results with NS-NC-34, for $\sqrt{v^\text{max}_l}=1.1$ p.u. $\forall l\in\mathcal{L}$ (left) and $\sqrt{v^\text{max}_l}=1.05$ p.u.  $\forall l\in\mathcal{L}$ (right). The load flow based on the dispatch from the R-OPF solution violates voltage bounds by at least 1\e{-2} p.u. for hours represented in red. For hours represented in blue, voltage bounds are satisfied in the load flow solution.}
\label{fig:vdual}
\vspace{-12pt}
\end{figure}

In our simulations, we observe that binding voltage constraints make R-OPF inexact only when dual variables of upper voltage bounds reach a certain threshold: see Figure \ref{fig:vdual} which relates relaxed constraint residuals to dual variables of upper voltage bounds, for NS-NC-34. In those cases, and in those cases only, the optimal dispatch according to R-OPF leads to a load flow which violates voltage bounds. In line with the observations from \cite{bunaiyan16}, we find that the direct relaxation remains exact with the IEEE123 feeder under normal operating constraints, including with binding upper voltage bounds. In the test case NS-NC-123, dual variables of upper voltage bounds reach a maximum value of 0.13 p.u., in which case the largest cone constraint residual is 5.9\e{-4} p.u..

It is also interesting to observe that, in all simulated cases where the R-OPF solution remains exact, there are binding upper bounds for active and reactive power withdrawals. In both our IEEE34 and IEEE123 test cases, several buses do not have PV stations (3 buses in IEEE34 and 28 buses in IEEE123). Thus, for these buses, the active and reactive power withdrawals are equal to the loads connected to them (which are considered inflexible) and constraint \eqref{OPF9} writes as an equality constraint. As a result, the ``upper bounds'' for active and reactive power withdrawals at those buses are binding in all simulated hours, including those where R-OPF retains exactness.

These results show that for practical applications the conditions for exactness of R-OPF are substantially milder than those proposed in the literature \cite{low14partII, farivar11, gan12}, although it might be more difficult to extract them mathematically. 


\subsubsection{Augmented relaxation from Gan et al. \cite{gan15}}

\begin{figure}[]
\centering
\includegraphics[width=3.5in]{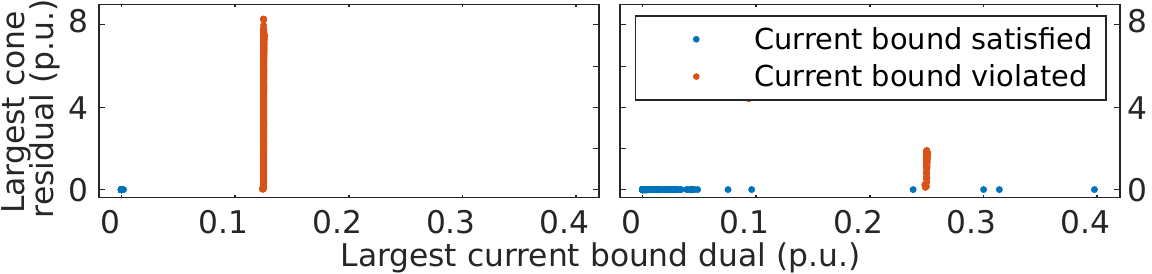}
\vspace{-15pt}
\caption{Inexactness of AR-OPF from Gan et al. \cite{gan15} as a function of dual variables of current bounds, with NS-C-123 (left, nominal current bound 4 p.u.) and NS-C-34 (right, nominal current bound 2 p.u.). The load flow based on the dispatch from the R-OPF solution violates current bounds by at least 1\e{-2} p.u. for hours represented in red. For hours represented in blue, both voltage and current bounds are satisfied in the load flow solution.}
\vspace{-12pt}
\label{fig:ldual}
\end{figure}

The authors in \cite{gan15} propose an augmented relaxation which ignores line shunts and current bounds (see Table \ref{tab:methods}), and provide a sufficient ex-ante condition for exactness based on network parameters. For the test cases NS-NC-34 and NS-NC-123, we observe that this condition is always satisfied, and that their AR-OPF retains exactness in all hours. These results confirm that the method from Gan et al. \cite{gan15} is applicable to networks where line shunts and current bounds are omitted. When including current bounds in the problem, we observe that the AR-OPF from \cite{gan15} can become inexact when current bounds are binding. In NS-C-123, current limits are binding in 19\% of hours, in which cases the relaxation is inexact and the optimal withdrawal setpoints from the AR-OPF lead to infeasible load flows (see Figure \ref{fig:ldual}, left). In NS-C-34, binding current bounds are however not necessarily leading to inexact solutions (see Figure \ref{fig:ldual}, right -- this case is simulated here with a current bound of 2 p.u. as current bounds of 4 p.u. are only rarely binding in NS-C-34). 
With both test feeders, we further observe that in all hours where the AR-OPF is inexact due to binding current bounds, the relaxation gap is zero. This means that the AR-OPF finds the globally optimal objective value, but the determined optimal setpoint is outside the feasible space. In Section \ref{res:4}, we show how we can render the AR-OPF from \cite{gan15} exact in all hours with binding current bounds by penalising current magnitudes in the objective function.
\subsubsection{Augmented relaxation from Huang et al. \cite{huang17}}
The authors in \cite{huang17} propose an augmented relaxation which readily includes current limits. Their definition of exactness differs slightly from \cite{gan15} and \cite{nick17}, in that they consider the AR-OPF exact when the optimal withdrawal setpoints for non-root nodes is a feasible dispatch for the network, even if the relaxed constraints \eqref{eq:SOC} are not binding. The authors prove that this property is verified when the active power import cost function is non-decreasing, which is a milder condition than that of Gan et al. \cite{gan15}. No ex-ante conditions on the network parameters are required, which is a benefit of the method. In cases where the A-OPF from Huang et al. \cite{huang17} is feasible (i.e., when there exist a solution with mild, non-simultaneous active and reactive reverse power flows on non-leaf lines, see Section \ref{res:1}), our simulations confirm that their AR-OPF provides feasible withdrawal setpoints for non-root nodes, even with binding current bounds.

\subsubsection{Augmented relaxation from Nick et al. \cite{nick17}}
The authors in \cite{nick17} propose an augmented relaxation which includes current limits, without constraints on reverse power flows. The authors provide sufficient ex-ante conditions for exactness based on network parameters. On test cases S-C-34 and S-C-123 (which include both line shunts and current bounds), we find that the AR-OPF from \cite{nick17} is always exact, but we observe that the ex-ante conditions for exactness are not satisfied. With S-C-34 we find that only condition C1 (see \cite{nick17}, p.6) is verified, while other conditions are not satisfied in any hour.\footnote{The authors in \cite{nick17} find that these ex-ante conditions are all satisfied for the IEEE34 feeder with full loading. We reproduce this result when omitting the low-voltage part of the IEEE34 feeder and assuming that the DERs do not have reactive power capability.} With S-C-123 we find that no condition is valid in any hour. This suggests that the method from Nick et al. \cite{nick17} has a broader range of application than indicated by the set of sufficient conditions provided by the authors.
\subsection{Application to reactive power control}
\label{res:3}

The solution methods from \cite{gan15,huang17,nick17} are applied by their authors to an OPF problem where the reactive power import $Q^t_1$ is neither constrained nor taken into account in the objective function. The OPF problem for radial networks however has applications for the optimal control of reactive power import/export of distribution grids. 
In order to test the applicability of these augmented SOCP relaxations to reactive power control, in this Section we test these with two different extensions:\begin{enumerate}
\item The reactive power import $Q^t_1$ is defined as a fixed parameter, positive for exporting inductive power, negative for importing inductive power, null for imposing reactive power balance at the root node. We let the objective function minimise the import of active power.
\item The objective function is the squared deviation of reactive power import/export from a target $Q^\text{ref}$, i.e. $(Q^t_1 - Q^\text{ref})^2$, with $Q^\text{ref}\in\mathbb{R}$.
\end{enumerate}
We test each method on test cases in line with the assumptions made by their authors. The method from \cite{huang17} is tested with variable capacitor bank capacitance to obtain feasible solutions on both test feeders. Results are summarised in Table \ref{tab:VC}, as the proportion of hours where an AR-OPF is inexact \textit{among the hours} where its respective A-OPF is feasible.

\begin{table}
\caption{Proportion of hours where AR-OPF is inexact among hours where A-OPF is feasible.}
\vspace{-10pt}
\label{tab:VC}
\begin{center}
\begin{tabular}{|c|c|c|c|c|c|c|}
\hline test grid & \multicolumn{3}{c|}{IEEE34} & \multicolumn{3}{c|}{IEEE123} \\
\hline \hspace{-5pt}configuration\hspace{-5pt} & NS-NC & NS-C$^\ast$ & S-C & NS-NC & NS-C$^\ast$ & S-C \\
\hline method & \cite{gan15} & \cite{huang17} & \cite{nick17} & \cite{gan15} & \cite{huang17} & \cite{nick17} \\
\hline $Q^t_1=0$ & 0.4\% & 0\% & 0.4\% & 0\% & 0\% & 0\%\\
\hline \hspace{-5pt}$Q^\text{ref}\hspace{-2pt}<\hspace{-2pt}\min\{Q^t_1\}$\hspace{-5pt} & 0\% & 0\% & 0\% & 0\% & 0\% & 0\% \\
\hline \hspace{-5pt}$Q^\text{ref}\hspace{-2pt}>\hspace{-2pt}\max\{Q^t_1\}$\hspace{-5pt} & 100\% & 100\% & 100\% & 100\% & 43\% & 100\% \\
\hline $Q^\text{ref}=0$ & 100\% & 100\% & 100\% & 100\% & 100\% & 100\% \\

\hline
\end{tabular}
\end{center}
\vspace{-3pt}
$^\ast$ Capacitor banks with variable capacitance.
\vspace{-10pt}
\end{table}

\subsubsection{Fixed reactive power import} We find that, when a feasible setpoint for the import of reactive power is enforced as a parameter, the augmented relaxations are inexact only in a minority of cases (see Table \ref{tab:VC}, first row, with a reactive power import enforced to zero). The relaxations from \cite{gan15} and \cite{nick17} are inexact on the IEEE34 test feeder only at times where the available active power from PV stations is particularly low and the load demand is high, representing less than 1\% of hours. All three relaxations are exact on the IEEE123 feeder. These results suggest that although exactness is not guaranteed for the augmented relaxations from \cite{gan15} and \cite{nick17} when the reactive power import is set to a fixed value, there might exist relatively mild conditions for exactness under such conditions. Exactness holds for \cite{huang17} in the two cases we studied.

\subsubsection{Minimising deviation from a target}
When setting a target for reactive power import in the problem's objective function, our simulation results differ depending on whether the target can be reached (i.e. there exists a feasible solution where $Q^t_1 = Q^\text{ref}$) or not.

When the target cannot be reached, our simulations show that the augmented relaxations from \cite{gan15,huang17,nick17} are:
\begin{itemize}
    \item exact, when the target is \textit{too low} to be reached, i.e. $Q^\text{ref}<\min\{Q^t_1\}$;
    \item inexact, when the target is \textit{too high} to be reached, i.e. $Q^\text{ref}>\max\{Q^t_1\}$.
\end{itemize}
When the target is too low, the solver has incentives to minimise current variables in the system, rendering constraint \eqref{eq:SOC} binding and leading to exact solutions. When the target is too high to be met, the solver has incentives to set current variables to higher values than those satisfying constraint \eqref{OPF4} to virtually increase the inductive effect of lines, which leads to inexact solutions. Results are summarised in Table \ref{tab:VC}, rows 2 and 3.

When the target for reactive power import can be met with a feasible solution, the augmented relaxations are inexact (see Table \ref{tab:VC}, row 4, for the case where $Q^\text{ref}=0$, which is a reachable target in these simulations) but the relaxation gap is zero. This means that there exist multiple optimal solutions to the relaxed problem, and at least one of those is feasible to the non-relaxed problem. Here is an example of how this can happen:\begin{itemize}
    \item let $\xi^\star\hspace{-2pt}=\hspace{-2pt}\{s^\star_l, v^\star_l, f^\star_l, S^{b\star}_l, S^{t\star}_l\}$ be an \textit{exact} optimal solution to the relaxed problem; the objective value $(Q^{t\star}_1 - Q^\text{ref})^2$ in this case is 0;
    \item for a given leaf line $L$, assume that there exists $\varepsilon$ so that $s'_L=s^\star_L+(j|z_L|^2 b_L-z_l)\varepsilon\in\mathcal{S}_L$ and $v'_L=v^\star_L+|z_l|^2\varepsilon\leq v^{\text{max}}_L$; let $f'_L=f^\star_L+\varepsilon$;
    \item the relaxed problem has a feasible solution $\xi'$ where $\{s_L,v_L,f_L\}=\{s'_L,v'_L,f'_L\}$ and all other variables are as in $\xi^\star$; with the solution $\xi'$, the objective value is also equal to 0, but constraint \eqref{eq:SOC} is not binding for line $L$, i.e. the solution is optimal but inexact, and the relaxation gap is zero.
\end{itemize}
In Section \ref{res:4}, we recover exactness in these cases by penalising current magnitudes in the objective function.


\subsection{Tightening relaxations with a zero relaxation gap}\label{res:4}
\begin{table}
\caption{Proportion of hours where relaxations are inexact among hours where: (row 1) the problem is feasible; (rows 2-4) the reactive power target can be met.}
\vspace{-10pt}
\label{tab:VD}
\begin{center}
\begin{tabular}{|c|c|c|c|c|c|}
\hline \multirow{2}{*}{method} & \multirow{2}{*}{test case} & \multicolumn{2}{c|}{IEEE34} & \multicolumn{2}{c|}{IEEE123} \\
\cline{3-6} & & $\epsilon=0$ & $\epsilon=0.01$ & $\epsilon=0$ & $\epsilon=0.01$ \\
\hline \cite{gan15} & NS-C$^\text{(a)}$ & 25\% & 0\% & 19\% & 0\% \\
\hline \cite{gan15} & NS-NC$^\text{(b)}$ & 100\% & 0\% & 100\% & 0\% \\
\hline \cite{nick17} & NS-C$^\text{(b)(c)}$ & 100\% & 0\% & 100\% & 0\% \\
\hline \cite{nick17} & S-C$^\text{(b)}$ & 100\% & 0.3\%$^\text{(d)}$ & 100\% & 0\% \\
\hline
\end{tabular}
\end{center}\vspace{-3pt}
\begin{itemize}
\item[(a)] OF$=P_1^t+\epsilon\sum_{l\in\mathcal{L}}f_l$, current limit 2 p.u.
\item[(b)] OF$={(Q_1^t-0)}^2+\epsilon\sum_{l\in\mathcal{L}}f_l$.
\item[(c)] Capacitor banks with variable capacitance.
\item[(d)] 0\% for $\epsilon=0.02$.
\end{itemize}
\vspace{-5pt}
\end{table}
\begin{figure}[]
\centering
\includegraphics[width=3.5in]{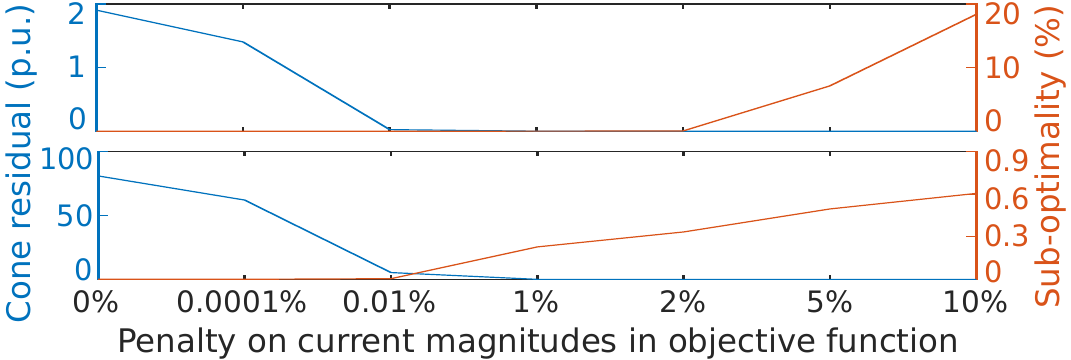}
\vspace{-15pt}
\caption{Relaxation tightness (measured as largest relaxed constraint residual throughout the year) and sub-optimality (measured as average increase in true objective throughout the year when introducing the penalty) as a function of a penalty on squared current magnitudes in the objective function. Above: AR-OPF from \cite{gan15} for NS-C-34 minimise active power import ($I_l^\text{max}=2$p.u.). Below: AR-OPF from \cite{nick17} for S-C-34 minimising reactive power import deviation from a target $Q^\text{ref}=0$.}
\label{fig:tight}
\vspace{-10pt}
\end{figure}
We observed in Section \ref{res:2} that the AR-OPF from Gan et al.~\cite{gan15} can lose exactness when current bounds are introduced in the problem, but that the relaxation gap is still zero in those cases. Similarly, we observed in \ref{res:3} that, when the problem's objective is to minimise deviation from a reachable target, all three augmented relaxations lose exactness but retain a zero relaxation gap. 
This means that in both situations, the AR-OPF has multiple optimal solutions, but the solver is not ``incentivised'' to choose an exact one among those. We investigate whether introducing a small penalty on current magnitudes in the objective function can incentivise the solver to minimise the residuals of the relaxed constraints and, thus, yield an exact solution. 
The objective function is augmented with a term $\epsilon\sum_{l\in\mathcal{L}}f_l$, with $\epsilon$ a chosen penalty factor. Table \ref{tab:VD} shows that a penalty of 1-2\% can render the augmented relaxations exact in the above-mentioned cases. Figure \ref{fig:tight} shows, for two examples, how the value of the penalty affects tightness and optimality. Tightness of the relaxations increases when the penalty increases. With a penalty of 1\%, the sub-optimality introduced in these examples is negligible (Figure \ref{fig:tight}, above) or small (0.3\%, Figure \ref{fig:tight}, below).

\section{Conclusion} \label{sec:ccl}
Convex relaxations of the Optimal Power Flow problem are essential not only for identifying the globally optimal solution (or a lower bound of it, if inexact) but also for enabling the use of OPF formulations in Bilevel Programming and Mathematical Programs with Equilibrium Constraints (MPEC), which are required for solving problems such as the TSO/DSO coordination or optimal network investment. This paper focuses on active distribution networks, and (i) introduces a \emph{unified} mathematical and simulation framework to assess the three most promising formulations for practical applications, (ii) it extends them to retain exactness and (iii) to consider reactive power injections. All three approaches \cite{gan15,huang17,nick17} consider a number of the practical limitations existing in real distribution networks, such as voltage limits and/or current limits, and supply \emph{sufficient} conditions under which they can guarantee that the determined optimal point is feasible to the original OPF problem. Our work is the first, to our knowledge, to compare and examine the extent of the applicability of the proposed methods under realistic conditions. We provide the simulation data, with the required modifications that facilitated the comparison, in an on-line appendix for future use by the interested reader \cite{appendix}. 

%
%

Our results show that the augmentations can lead to a non-negligible shrinkage of the feasible space of the original OPF problem and can incur sub-optimality and infeasibility under realistic operating conditions. Despite that, we show that the methods do retain exactness in cases that are not encompassed by their initially proposed sufficient conditions. For inexact cases with a zero relaxation gap, we include a small penalty on current magnitudes in the objective function that substantially increases the range of cases where an exact solution is recovered. Finally, we specify the use of a numerical threshold to assess exactness: in our simulation results, it appears that when all residuals of the relaxed constraints are below $1 \times 10^{-2}$ p.u., the relaxation is exact. 

From the three methods we examined in detail, we found that the solution from Nick et al. \cite{nick17} is the most promising for general radical networks. If line shunts are negligible, then the solution method from Gan et al. \cite{gan15} holds potential, as it is simpler. While this method does not originally include current constraints, our results suggest that exactness holds if we add the mild penalty on current magnitudes we have introduced. Finally, if the network experiences only mild reverse flows (which does not happen often in active distribution networks) and has negligible line shunts, the solution method from Huang et al. \cite{huang17} could also be applied, as it can guarantee exactness under milder sufficient conditions.

Future research directions will focus on the integration of the method from Nick et al. \cite{nick17} to bilevel programs for strategic bidding and TSO/DSO coordination. Promising directions also include (i) the formulation and proof of milder ex-ante conditions for exactness of the method from Nick et al. \cite{nick17}, with applications to both active and reactive power control, (ii) conditions for exactness of the method from Gan et al. \cite{gan15} under binding current constraints, so as to better represent the range of networks and conditions these methods can be applied to, and (iii) multi-phase distribution networks \cite{zhou2019sufficient, zhao2017convex}.

\bibliographystyle{IEEEtran}
\bibliography{biblio}
\end{document}